\documentclass{article}
\usepackage[utf8]{inputenc}
\usepackage{graphicx}
\usepackage{amsmath}
\usepackage{amssymb}
\usepackage{enumerate}
\usepackage{latexsym}
\usepackage{epsfig}
\usepackage{amsthm}
\usepackage{authblk}
\usepackage{multirow}

\newtheorem{definition}{Definition}

\newtheorem{theorem}{Theorem}

\title{
On the linear space
of the two-sided generalized Fibonacci sequences} 

\author{
Martin Bunder \\
School of Mathematics and Applied Statistics,\\
University of Wollongong, 
Australia, \\ {\tt{martin.bunder@uow.edu.au}}
\and
Joseph Tonien \\
School of Computing and Information Technology,\\
University of Wollongong,
Australia, \\ {\tt{joseph.tonien@uow.edu.au}}
}

\date{}
    
\begin{document}

\maketitle

\begin{abstract}
In this paper, we study the
linear space of all two-sided generalized Fibonacci sequences $\{F_n\}_{n \in \mathbb{Z}}$  that satisfy the recurrence equation of order $k$:
$F_n = F_{n-1} + F_{n-2} + \dots + F_{n-k}$.
We give two types of explicit formula, one is based on  generalized binomial coefficients and the other based on  generalized multinomial coefficients.
\end{abstract}

{\em AMS Classification Numbers:}
11B37, 11B39, 47B37


\section{Introduction}

The Fibonacci sequence,
$F_0=0$, $F_1 = 1$,
$F_n = F_{n-1} + F_{n-2}$,
have been generalized in many ways.
One of the generalizations~\cite{Kolodner1965,Cusick1968,Melham2004} is
to change the recurrence equation to
$F_n = \alpha F_{n-1} + \beta F_{n-2}$,
thus keeping the characteristic equation
remained in order 2.
Another common generalization
is to extend the recurrence equation
to a higher order.
For a fixed integer $k \geq 2$,
a sequence is called a Fibonacci sequence of order $k$ if it satisfies the following recurrence equation
\begin{equation}
\label{GF:eqn1}
    F_n = F_{n-1} + F_{n-2} + \dots + F_{n-k}.
\end{equation}
For some particular values of $k$,
the sequence has a special name.
It is called a
tribonacci sequence, 
a tetranacci sequence
and a pentanacci sequence
for $k=3, 4, 5$, respectively.

A Fibonacci sequence of order $k$ is uniquely determined by a list of values of $k$ consecutive terms. For instance, 
if the values of $F_0, F_1, \dots, F_{k-1}$ are given then using the recurrence equation~(\ref{GF:eqn1}), we can work out the values of all other terms $F_n$ for $n \geq k$, as well as for negative indices $n<0$. Here is an example of a Fibonacci sequence of order 5:
\begin{align*}
\dots,
F_{-4} = -2,
F_{-3} = 7,
F_{-2} = -3,
F_{-1} = -4,
\\
F_0 = {\bf 3},
F_1 = {\bf 1},
F_2 = {\bf 4},
F_3 = {\bf 1},
F_4 = {\bf 5},
F_5 = 14,
F_6 = 25, \dots
\quad
.
\end{align*}

Since we have $F_0=0$ and $F_1=1$ in the original Fibonacci sequence, there are two common ways to set the initial conditions:
(i) $F_0=F_1=\dots=F_{k-2}=0$, $F_{k-1}=1$ as in~\cite{Miles1960,Gabai1970,Philippou1982,Lee2001,Chaves2014,Ddamulira2018};
or (ii) $F_0=0$, $F_1=\dots=F_{k-2}=F_{k-1}=1$ as in~\cite{Lengyel_Marques_2017,Sobolewski2017,Bunder_Tonien2020}. Another initial condition 
$F_0=F_1=\dots=F_{k-1}=1$ appears in Ferguson~\cite{Ferguson1966} arisen in the study of polyphase merge-sorting. Various formulas have been found 
for Fibonacci sequences with these three initial conditions
which can be grouped into three types:
Binet formula~\cite{Dence1987,Lee2001},
binomial coefficients~\cite{Ferguson1966,Benjamin2014}
and multinomial coefficents~\cite{Miles1960,Lee2001}.
We note that these formulas of $F_n$ are only restricted to the integer indices $n \geq 0$. 
The Binet  type of formula is algebraic in nature and remains valid when we extend to negative indices $n<0$.
However, formulas involved binomial coefficients 
and multinomial coefficents are limited to non-negative indices and it is not trivial to extend to negative indices.

While most authors only consider sequences $F_n$ with $n \geq 0$,
in this paper, we will study two-sided sequences. Those are sequences $\{F_n\}$
where the index $n \in \mathbb{Z}$,
that is, we allow $n$ to be a {\em negative integer}.
Instead of looking for explicit formula
for a Fibonacci sequence with a particular initial condition, our aim is to 
find explicit formulas for
a general Fibonacci sequence
that has an {\em arbitrary 
initial condition} $(F_0,F_1,\dots,F_{k-1})$.
To do that, we consider
the set of all Fibonacci sequences of order $k$.
This
forms a $k$-dimensional linear space.
We will study the standard basis of this linear space
which is denoted by
$B^{(0)}, B^{(1)}, \dots, B^{(k-1)}$.
For $0 \leq j \leq k-1$, each $B^{(j)}$ is a Fibonacci sequence
whose initial values are all zero except
$B^{(j)}_j=1$.
We will find explicit formula for the basis sequences $B^{(0)}, B^{(1)}, \dots, B^{(k-1)}$, and thus, 
any Fibonacci sequence $F$ can be determined by a linear combination $F=F_0 B^{(0)} + F_1 B^{(1)}  +  \dots + F_{k-1} B^{(k-1)}$.

Our aim is to find explicit formulas
for two-sided Fibonacci sequences
that are expressed in terms of binomial coefficients
and multinomial coefficients, respectively.
Since the classical binomial coefficients
and multinomial coefficients are only associated with non-negative integers,
to use these for our two-sided sequences we need to extend the binomial notation
and multinomial notation to include negative integers.
To this end, 
we extend the binomial notation ${n \choose i}$ to negative values of $n$ and $i$, writing this as $\left\langle {n \choose i} \right\rangle$. 
Subjected to the two conditions
$\left\langle {n \choose n} \right\rangle=1$
and
$\left\langle {{n-1} \choose {i}} \right\rangle 
+\left\langle {{n-1} \choose {i-1}} \right\rangle = \left\langle {{n} \choose {i}} \right\rangle$, the latter is called {\em the Pascal Recursion equation},
the value of the generalized binomial notation is uniquely determined.
In Theorem~\ref{Bkj:theorem2}, we will show that
\begin{align*}
B_n^{(j)} = & - \sum_{i \in \mathbb{Z}}{ (-1)^i \left\langle {{n-ik} \choose {i-1}} \right\rangle 2^{n+1-i(k+1)}}
\\
&
+ \sum_{i \in \mathbb{Z}}{ (-1)^i  \left\langle {{n-j-1-ik} \choose {i-1}} \right\rangle 2^{n-j-i(k+1)}}
\mbox{ for all } n \in \mathbb{Z}.
\end{align*}

We extend the multinomial notation  ${n \choose {i_1, i_2, \dots, i_t}}$ to negative values of $n$ and $i_1,\dots,i_t$, writing this as $\left\langle {n \choose {i_1, i_2, \dots, i_t}} \right\rangle$. 
The generalization is done as follows.

Using the generalized binomial notation we extend the traditional multinomial notation 
$$
{n \choose {i_1, i_2, \dots, i_k}}=
{n \choose {i_2 + \dots + i_{t}}}
{{i_2 + \dots + i_{t}} \choose {i_3 + \dots + i_{t}}}
\dots 
{{i_{t-2} + i_{t-1} + i_{t}} \choose {i_{t-1} + i_t}}
{{i_{t-1} + i_{t}} \choose {i_t}},
$$
to
\begin{align*}
&\left\langle {{n} \choose {i_1, i_2, \dots, i_t}} \right\rangle
=
\left\langle {{n} \choose {i_2 + \dots + i_{t}}} \right\rangle
\left\langle {{i_2 + \dots + i_{t}} \choose {i_3 + \dots + i_{t}}} \right\rangle
\dots 
\left\langle {{i_{t-2} + i_{t-1} + i_{t}} \choose {i_{t-1} + i_t}} \right\rangle
\left\langle {{i_{t-1} + i_{t}} \choose {i_t}} \right\rangle
.
\end{align*}

Using this generalized multinomial notation,
in Theorem~\ref{MilesExtend}, we will show that
$$
B_{n}^{(j)}
=
\sum_{n-k-j \leq  a_1 + 2 a_2 + \dots + k a_k \leq n-k}{\left\langle {{a_1+a_2 + \dots + a_k} \choose {a_1, a_2, \dots, a_k}} \right\rangle}, \mbox{ for all } n \in \mathbb{Z}.
$$

The rest of the paper is organised as follows.
In section~\ref{GF:sec2}, we study the linear space of Fibonacci sequences of order $k$ in general, especially looking at the linear automorphisms of this space. Formulas based on the generalized binomial notation are derived in section~\ref{GF:sec3}.
Formulas based on the generalized multinomial notation are derived in section~\ref{GF:sec4}.
Finally, in section~\ref{GF:sec5}, we remark on how the generalized Fibonacci sequences are related to a tiling problem.

\section{The Fibonacci linear space of order $k$}
\label{GF:sec2}

\begin{definition}
Let $k \geq 2$ be a fixed integer.
A sequence $\{F_n\}_{n \in {\mathbb{Z}}}$
is called a Fibonacci sequence of order $k$
if it satisfies
 the following recurrence equation
\begin{equation}
\label{kFiboRecurrence}
    F_n = F_{n-1} + F_{n-2} + \dots + F_{n-k}, \mbox{ for all } n \in \mathbb{Z}.
\end{equation}
\end{definition}

We can see that, given $k$ values $(F_0, F_1, \dots, F_{k-1})$,
then using the Fibonacci recurrence equation~(\ref{kFiboRecurrence}),
all other values $F_n$ for $n \in \mathbb{Z}$ are determined uniquely. We will refer to $(F_0, F_1, \dots, F_{k-1})$
as the initial values of the sequence.
The set of all Fibonacci sequences of order $k$ forms a $k$-dimensional vector space
(either over the field $\mathbb{R}$ or $\mathbb{C}$).
We will use $\mathsf{Fibonacci}^{(k)}$ to denote
this vector space of all Fibonacci sequences of order $k$.
We now define the standard basis 
for the Fibonacci vector space $\mathsf{Fibonacci}^{(k)}$.

\begin{definition}
Let $k \geq 2$ be a fixed integer.
For each integer $0 \leq j \leq k-1$,
the sequence $B^{(j)} \in \mathsf{Fibonacci}^{(k)}$
is defined by the initial values
$$
B^{(j)}_n =
\begin{cases}
			0, & \mbox{ if } 0 \leq n \leq k-1 \mbox{ and } n \neq j\\
            1, &
            \mbox{ if } n = j.
		 \end{cases}
$$
\end{definition}

The special sequences $B^{(0)}, B^{(1)}, \dots, B^{(k-1)}$ defined above form a standard basis for the space $\mathsf{Fibonacci}^{(k)}$.
Any member of this Fibonacci vector space is a linear combination of the standard basis and we have the following theorem.

\begin{theorem}
\label{Bkj:theorem0}
Let $k \geq 2$ be a fixed integer.
Let $\{F_n\}_{n \in {\mathbb{Z}}}$ be 
a Fibonacci sequence of order $k$.
Then
$$
F_n = \sum_{j=0}^{k-1} {B^{(j)}_n F_j}
\mbox{ for all } n \in \mathbb{Z}.
$$
\end{theorem}

By Theorem~\ref{Bkj:theorem0},
we can see that in order to determine an explicit formula for any
Fibonacci sequence $\{F_n\}_{n \in {\mathbb{Z}}}$,
it suffices to derive formula for 
the $k$ basis sequences $B^{(0)}, B^{(1)}, \dots, B^{(k-1)}$.

\subsection{Linear operators on the Fibonacci space}

Here we list some standard linear operators 
on two-sided sequences.

\begin{itemize}
    \item Identity operator
    $\mathtt{I}$.
    
    \item Left shift operator $\mathtt{L}$: $\mathtt{L}(X)=Y$ iff $Y_n = X_{n+1}$
    for all $n \in {\mathbb{Z}}$.
    
    \item Right shift operator $\mathtt{R}$:
    $\mathtt{R}(X)=Y$ iff $Y_n = X_{n-1}$
    for all $n \in {\mathbb{Z}}$.
    The left shift and the right shift are inverse of each other: 
    $\mathtt{L} \mathtt{R} = \mathtt{R} \mathtt{L} = \mathtt{I}$.
    
    \item Forward difference operator $\Delta$:
    $\Delta(X)=Y$ iff $Y_n = X_{n+1} - X_n$
    for all $n \in {\mathbb{Z}}$.
    Here $\Delta = \mathtt{L} - \mathtt{I}$.
    
    \item Backward difference operator $\nabla$: $\nabla(X)=Y$ iff $Y_n = X_{n} - X_{n-1}$
    for all $n \in {\mathbb{Z}}$.
    Here $\nabla = \mathtt{I} - \mathtt{R} = \mathtt{I} - \mathtt{L}^{-1}$, $\mathtt{L} \nabla = \Delta$
    and $\mathtt{R} \Delta = \nabla$.
\end{itemize}
 
We have the following theorem
concerning the above operators.
 
 \begin{theorem}
 \label{thm_operator}
 All operators
 $\mathtt{I}$, $\mathtt{L}$,
 $\mathtt{R}$, $\Delta$ and
$\nabla$ when restricted to the space $\mathsf{Fibonacci}^{(k)}$ are linear automorphisms $\mathsf{Fibonacci}^{(k)} \to \mathsf{Fibonacci}^{(k)}$
and satisfy the following relations:
 
(i) $\mathtt{L}^{k} = \mathtt{I} + \mathtt{L} + \mathtt{L}^2 + \dots + \mathtt{L}^{k-1}
.
$
 
(ii) $\mathtt{R} = \mathtt{L}^{-1} = - \mathtt{I} - \mathtt{L} - \mathtt{L}^2 - \dots 
 - \mathtt{L}^{k-2} + \mathtt{L}^{k-1}$
 
 (iii)
 $\mathtt{R}^{k} = \mathtt{I} - \mathtt{R} - \mathtt{R}^{2} - \dots - \mathtt{R}^{k-1}
 .
 $
 
 (iv)
 $\mathtt{L} = \mathtt{R}^{-1} = \mathtt{I} + \mathtt{R} + \mathtt{R}^{2} + \dots + \mathtt{R}^{k-1}
 .
 $

(v)
 $\mathtt{L}^{k+1}= 2 \mathtt{L}^{k} - \mathtt{I}
 .
 $
 
 (vi)
 $\mathtt{R}^{k+1} = 2 \mathtt{R} - \mathtt{I}
 .
 $
 
 (vii)
 $\Delta
     ( \mathtt{I} + (k-1) \mathtt{R} + (k-2) \mathtt{R}^2
     + (k-3) \mathtt{R}^3 + \dots + 
      2 \mathtt{R}^{k-2} +  \mathtt{R}^{k-1})
     = (k-1) \mathtt{I}
     .
     $

 (viii)
$\nabla
     (k \mathtt{I} + (k-1) \mathtt{R}
     + (k-2) \mathtt{R}^2 + \dots + 
      2 \mathtt{R}^{k-2} + \mathtt{R}^{k-1})
= (k-1) \mathtt{I}
.
$
 
  (ix)
 $
    \sum_{i=0}^{k}
    {{k+1} \choose {i+1}}
    \frac{k-1-2i}{k+1}
    \Delta^{i}
    = 0
    .
 $

 (x) $
   (k-1)\mathtt{I}
     + \sum_{i=1}^{k}  {{k+1} \choose {i+1}}  (-1)^{i} \nabla^i
     =0
     .
 $

 \end{theorem}
 {\em Proof.}
 It is easy to see that all these operators $\mathtt{I}$, $\mathtt{L}$,
 $\mathtt{R}$, $\Delta$ and
$\nabla$ are linear. Each maps a Fibonacci sequence to another Fibonacci sequence. The bijectivity of  
 $\mathtt{I}$, $\mathtt{L}$,
 $\mathtt{R}$ is obvious, whereas, the bijectivity of $\Delta$ and
$\nabla$ follows from (vii) and (viii), respectively.

(i) For any 
$X \in \mathsf{Fibonacci}^{(k)}$,
let $(\mathtt{I} + \mathtt{L} + \mathtt{L}^2 + \dots + \mathtt{L}^{k-1})(X) = Y$
then $Y_n = X_n + X_{n+1} + X_{n+2} + \dots + X_{n+k-1} = X_{n+k}$,
therefore,
$Y = \mathtt{L}^{k}(X)$.
This proves that, restricted to the linear space $\mathsf{Fibonacci}^{(k)}$,
$\mathtt{I} + \mathtt{L} + \mathtt{L}^2 + \dots + \mathtt{L}^{k-1} = \mathtt{L}^{k}$.

(ii) For any 
$X \in \mathsf{Fibonacci}^{(k)}$,
let $(- \mathtt{I} - \mathtt{L} - \mathtt{L}^2 - \dots 
 - \mathtt{L}^{k-2} + \mathtt{L}^{k-1})(X)=Y$
 then $Y_n = -X_n - X_{n+1} - X_{n+2} - \dots - X_{n+k-2} + X_{n+k-1} = X_{n-1}$.
Hence, $Y = \mathtt{R}(X)$,
 and therefore,
$- \mathtt{I} - \mathtt{L} - \mathtt{L}^2 - \dots 
 - \mathtt{L}^{k-2} + \mathtt{L}^{k-1} = \mathtt{R} =  \mathtt{L}^{-1}$.
 
 (iii)
 For any 
$X \in \mathsf{Fibonacci}^{(k)}$,
let $(\mathtt{I} - \mathtt{R} - \mathtt{R}^{2} - \dots - \mathtt{R}^{k-1})(X)=Y$
then $Y_n = X_n - X_{n-1} - X_{n-2} - \dots - X_{n-k+1} = X_{n-k}$.
Hence, $Y = \mathtt{R}^{k}(X)$, and therefore,
 $\mathtt{I} - \mathtt{R} - \mathtt{R}^{2} - \dots - \mathtt{R}^{k-1}= \mathtt{R}^{k}.$
 
 (iv)
 For any 
$X \in \mathsf{Fibonacci}^{(k)}$,
let $(\mathtt{I} + \mathtt{R} + \mathtt{R}^{2} + \dots + \mathtt{R}^{k-1})(X)=Y$ then
$Y_n = X_n + X_{n-1} + X_{n-2} + \dots + X_{n-k+1} = X_{n+1}$.
Hence, $Y = \mathtt{L}(X)$, and
therefore,
 $\mathtt{I} + \mathtt{R} + \mathtt{R}^{2} + \dots + \mathtt{R}^{k-1}=\mathtt{L} = \mathtt{R}^{-1}$.

 (v)
 By (i),
 $\mathtt{L}^{k+1} = \mathtt{L} \, \mathtt{L}^{k}= \mathtt{L} (\mathtt{I} + \mathtt{L} + \mathtt{L}^2 + \dots + \mathtt{L}^{k-1})= \mathtt{L}
 + \mathtt{L}^2 + \dots + \mathtt{L}^{k-1} 
 + \mathtt{L}^{k}= 
 (\mathtt{I} + \mathtt{L} + \mathtt{L}^2 + \dots + \mathtt{L}^{k-1}) + \mathtt{L}^{k} - \mathtt{I}= \mathtt{L}^{k} + \mathtt{L}^{k} - \mathtt{I}=
 2 \mathtt{L}^{k} - \mathtt{I}
 $.
 
 (vi)
 By (iii), $\mathtt{R}^{k+1}=
 \mathtt{R} \mathtt{R}^{k}=
 \mathtt{R} (\mathtt{I} - \mathtt{R} - \mathtt{R}^{2} - \dots - \mathtt{R}^{k-1})=
 \mathtt{R} - \mathtt{R}^2 - \mathtt{R}^3 - \dots - \mathtt{R}^{k-1} - \mathtt{R}^k =
 \mathtt{R} - \mathtt{R}^2 - \mathtt{R}^3 - \dots - \mathtt{R}^{k-1} -
 (\mathtt{I} - \mathtt{R} - \mathtt{R}^{2} - \dots - \mathtt{R}^{k-1})= 2 \mathtt{R} - \mathtt{I}$.
 
 (vii) We have
\begin{align*}
&
\Delta
     ( \mathtt{I} + (k-1) \mathtt{R} + (k-2) \mathtt{R}^2
     + (k-3) \mathtt{R}^3 + \dots + 
      2 \mathtt{R}^{k-2} +  \mathtt{R}^{k-1})
     \\
  & =  (\mathtt{L} - \mathtt{I})
     ( \mathtt{I} + (k-1) \mathtt{R} + (k-2) \mathtt{R}^2
     + (k-3) \mathtt{R}^3 + \dots + 
      2 \mathtt{R}^{k-2} +  \mathtt{R}^{k-1})
     \\
     &= \mathtt{L} + (k-2) 
     \mathtt{I}
     - \mathtt{R}
     -\mathtt{R}^2
     - \dots -
     \mathtt{R}^{k-2}
     - \mathtt{R}^{k-1}
     \\
     &= (k-1) \mathtt{I}
    \quad  \mbox{ by (iv).}
\end{align*}  

 (viii) We have
\begin{align*}
&
\nabla
     (k \mathtt{I} + (k-1) \mathtt{R}
     + (k-2) \mathtt{R}^2 + \dots + 
      2 \mathtt{R}^{k-2} + \mathtt{R}^{k-1})
     \\
  & =  (\mathtt{I} - \mathtt{R})
     (k \mathtt{I} + (k-1) \mathtt{R}
     + (k-2) \mathtt{R}^2 + \dots + 
     2 \mathtt{R}^{k-2} + \mathtt{R}^{k-1})
     \\
     &= k \mathtt{I} - 
     \mathtt{R}
     - \mathtt{R}^2 
     - \dots -
     \mathtt{R}^{k-1}
     - \mathtt{R}^k
     \\
     &= (k-1) \mathtt{I}
    \quad  \mbox{ by (iii).}
\end{align*}    
 
 (ix) 
 Substituting $\mathtt{L} = \mathtt{I} + \Delta$
 into (i), we have
 \begin{align*}
     (\mathtt{I} + \Delta)^{k} &= \mathtt{I} + (\mathtt{I} + \Delta) + (\mathtt{I} + \Delta)^2 + \dots + (\mathtt{I} + \Delta)^{k-1}
     \\
    \sum_{i=0}^{k} {k \choose i} \Delta^{i} &= 
    \sum_{j=0}^{k-1}
    \sum_{i=0}^{j} {j \choose i}
    \Delta^{i}=
    \sum_{i=0}^{k-1}
    \sum_{j=i}^{k-1} {j \choose i}
    \Delta^{i}
    =
    \sum_{i=0}^{k-1}
      {k \choose {i+1}}
    \Delta^{i}
    .
 \end{align*}
 Therefore,
 \begin{align*}
     \Delta^{k} &= 
    \sum_{i=0}^{k-1}
    \left(
    {k \choose {i+1}}
    -
    {k \choose i}
    \right)
    \Delta^{i}
    = 
    \sum_{i=0}^{k-1}
    {{k+1} \choose {i+1}}
    \frac{k-1-2i}{k+1}
    \Delta^{i}
    .
 \end{align*}
 
 (x)
 Substituting $\mathtt{R} = \mathtt{I} - \nabla$ into (iii), we have
 \begin{align*}
     (\mathtt{I} - \nabla)^k &= \mathtt{I} - (\mathtt{I} - \nabla) 
     - (\mathtt{I} - \nabla)^2 - \dots -
     (\mathtt{I} - \nabla)^{k-1}
     .
   \end{align*}  
   So
    \begin{align*}
     \sum_{i=1}^{k} {k \choose i} (-\nabla)^i
     &= -\sum_{j=1}^{k-1} \sum_{i=0}^{j} {j \choose i} (-\nabla)^i
     =-(k-1)\mathtt{I} - \sum_{i=1}^{k-1} \sum_{j=i}^{k-1} {j \choose i} (-\nabla)^i
     \\
     &
     =-(k-1)\mathtt{I} - \sum_{i=1}^{k-1} {k \choose {i+1}} (-\nabla)^i
     .
 \end{align*}
 Therefore,
 \begin{align*}
     (- \nabla)^k &
     =-(k-1)\mathtt{I} - \sum_{i=1}^{k-1} \left( {k \choose {i+1}} + {k \choose i} \right) (-\nabla)^i
     \\
     &
     =-(k-1)\mathtt{I} - \sum_{i=1}^{k-1}  {{k+1} \choose {i+1}}  (-\nabla)^i
 \end{align*}
 and
 \begin{align*}
    \sum_{i=1}^{k}  {{k+1} \choose {i+1}}  (-\nabla)^i
   &  =-(k-1)\mathtt{I}.
    \quad
 \blacksquare
 \end{align*}

 \begin{theorem}
 \label{BSthm}
 Denote $S = B^{(0)}+ B^{(1)}+ \dots+ B^{(k-1)}  \in \mathsf{Fibonacci}^{(k)}$.
 We have
 
 (i)
 $B^{(j)} -  B^{(j-1)}= 
\mathtt{R}^j(B^{(0)})
$
for all $1 \leq j \leq k-1$.

 (ii)
 $B^{(j)} = 
\sum_{i=0}^{j} \mathtt{R}^i(B^{(0)})
$
for all $0 \leq j \leq k-1$.

 (iii)
 $B^{(0)}= \mathtt{R}(B^{(k-1)})$
 and $B^{(k-1)}= \mathtt{L}(B^{(0)})$.
 
 (iv) $
B^{(j)}=\sum_{i=0}^{j} 
\mathtt{R}^{i+1}(B^{(k-1)})
$
for all $0 \leq j \leq k-1$.

(v)
$S = (k \, \mathtt{I} + (k-1) \, \mathtt{R}
+ (k-2) \, \mathtt{R}^2 + \dots + 
\mathtt{R}^{k-1})(B^{(0)})$.

(vi)
$\nabla(S) = (k-1) B^{(0)}$.

(vii)
$(\mathtt{I}- \mathtt{R}^{j+1})(S)
=(k-1) B^{(j)}$
for all $0 \leq j \leq k-1$.
 \end{theorem}
 {\em Proof.}
 (i)
 Both $B^{(j)} -  B^{(j-1)}$ and 
$\mathtt{R}^j(B^{(0)})$
 are members of $\mathsf{Fibonacci}^{(k)}$
 and their initial values are equal,
 therefore,
$B^{(j)} -  B^{(j-1)}= 
\mathtt{R}^j(B^{(0)})
$.
 
 (ii) 
 It follows from (i).

 (iii)
 By (ii),
 $B^{(k-1)} = 
\sum_{i=0}^{k-1} \mathtt{R}^i(B^{(0)})
$
and since $\mathtt{L} = \mathtt{R}^{-1}=
\mathtt{I} + \mathtt{R} + \mathtt{R}^2 + \dots + \mathtt{R}^{k-1}$ (Theorem~\ref{thm_operator}(iv)),
we have $B^{(k-1)} = 
\mathtt{L}(B^{(0)})
$
and so 
 $B^{(0)}= \mathtt{R}(B^{(k-1)})$.

  (iv) It follows from (ii) and (iii).
  
  (v) It follows from (ii).
  
  (vi) It follows from (v) and  Theorem~\ref{thm_operator}(viii).
  
  (vii)
  We have
  \begin{align*}
      (k-1) B^{(j)}
      & =(k-1) \sum_{i=0}^{j} \mathtt{R}^i(B^{(0)})
      \quad \mbox{ by (ii)}
      \\
     & = \sum_{i=0}^{j} \mathtt{R}^i(\nabla(S))
      \quad \mbox{ by (vi)}
       \\
     & = \sum_{i=0}^{j} (\mathtt{R}^i (1-\mathtt{R})) (S)
     =(1-\mathtt{R}^{j+1})(S).
  \end{align*}
  
Another direct way to prove (vii)
is by observing that 
both $(k-1) B^{(j)}$ and 
$(1-\mathtt{R}^{j+1})(S)$
 are members of $\mathsf{Fibonacci}^{(k)}$
 and their initial values are equal.
   \quad
$\blacksquare$ 

\section{Explicit formulas based on binomials}

\label{GF:sec3}

In this section, we will derive 
explicit formula
for 
the two-sided Fibonacci basis sequences
$B^{(0)}, B^{(1)}, \dots, B^{(k-1)}$
expressed in terms of binomial coefficients.
Since the traditional binomial notation is
associated with non-negative integers,
to use these for our two-sided sequences we need to extend the binomial notation
 to include negative integers.
To this end, 
we extend the binomial notation ${n \choose i}$ to negative values of $n$ and $i$.

The binomial notation ${n \choose i}$ can
be generalized to $\left\langle {{n} \choose {i}} \right\rangle$ for all integers $n$ and $i$ by enforcing two conditions:
\begin{itemize}
    \item $\left\langle {{n} \choose {n}} \right\rangle = 1$ for all $n \in \mathbb{Z}$; and
    
\item
Pascal Recursion relation
\begin{equation}
\label{genBinoEqn} 
\left\langle {{n-1} \choose {i}} \right\rangle 
+\left\langle {{n-1} \choose {i-1}} \right\rangle = \left\langle {{n} \choose {i}} \right\rangle .
\end{equation}
 \end{itemize}
 
 With these two conditions, 
 $\left\langle {{n} \choose {i}} \right\rangle$ is uniquely 
 determined as 
\begin{align}
\label{nchoosei1}
 \left\langle {{n} \choose {i}} \right\rangle
 & =
  \begin{cases}
  \frac{n^{\underline{n-i}}}{(n-i)!}=
\frac{n(n-1)(n-2)\dots(i+1)}{(n-i)!}, & \text{if } n \geq i\\
            0, & \text{otherwise}
		 \end{cases}
		 \\
\label{nchoosei2}
&= \begin{cases}
{n \choose i}, & \text{if } n \geq i \geq 0\\
(-1)^{i+n} {{-i-1} \choose {-n-1}}
, & \text{if } -1 \geq n \geq i\\
            0, & \text{otherwise}
		 \end{cases}
   .
\end{align}
 
 Refer to~\cite{Loeb1992, Loeb1992B} for detailed discussion on various generalizations of binomial notation.
 The following table shows some values of 
 $\left\langle {{n} \choose {i}} \right\rangle
 $:

\begin{tabular}{ |lr|rrrrrrrrrrrrr|  }
\hline
 \multicolumn{2}{|c|}{\multirow{2}{2em}{$\left\langle {{n} \choose {i}} \right\rangle$}} & 
 \multicolumn{13}{|c|}{$i$} \\
  & & $-6$ & $-5$ & $-4$ & $-3$ & $-2$ & $-1$ & $0$ & $1$ & $2$ & $3$ & $4$ & $5$ & $6$\\
 \hline
 \multirow{9}{1mm}{$n$}   
 & $6$ & $0$ & $0$ & $0$ & $0$ & $0$ & $0$ & $1$ & $6$ & $15$ & $20$ & $15$ & $6$ & $1$\\
 & $5$ & $0$ & $0$ & $0$ & $0$ & $0$ & $0$ & $1$ & $5$ & $10$ & $10$ & $5$ & $1$ & $0$\\
 & $4$ & $0$ & $0$ & $0$ & $0$ & $0$ & $0$ & $1$ & $4$ & $6$ & $4$ & $1$ & $0$ & $0$\\
 & $3$ & $0$ & $0$ & $0$ & $0$ & $0$ & $0$ & $1$ & $3$ & $3$ & $1$ & $0$ & $0$ & $0$\\
 & $2$ & $0$ & $0$ & $0$ & $0$ & $0$ & $0$ & $1$ & $2$ & $1$ & $0$ & $0$ & $0$ & $0$\\
 & $1$ & $0$ & $0$ & $0$ & $0$ & $0$ & $0$ & $1$ & $1$ & $0$ & $0$ & $0$ & $0$ & $0$\\
 & $0$ & $0$ & $0$ & $0$ & $0$ & $0$ & $0$ & $1$ & $0$ & $0$ & $0$ & $0$ & $0$ & $0$\\
 & $-1$ & $-1$ & $1$ & $-1$ & $1$ & $-1$ & $1$ & $0$ & $0$ & $0$ & $0$ & $0$ & $0$ & $0$\\
& $-2$ & $5$ & $-4$ & $3$ & $-2$ & $1$ & $0$ & $0$ & $0$ & $0$ & $0$ & $0$ & $0$ & $0$\\
& $-3$ & $-10$ & $6$ & $-3$ & $1$ & $0$ & $0$ & $0$ & $0$ & $0$ & $0$ & $0$ & $0$ & $0$\\
& $-4$ & $10$ & $-4$ & $1$ & $0$ & $0$ & $0$ & $0$ & $0$ & $0$ & $0$ & $0$ & $0$ & $0$\\
& $-5$ & $-5$ & $1$ & $0$ & $0$ & $0$ & $0$ & $0$ & $0$ & $0$ & $0$ & $0$ & $0$ & $0$\\
& $-6$ & $1$ & $0$ & $0$ & $0$ & $0$ & $0$ & $0$ & $0$ & $0$ & $0$ & $0$ & $0$ & $0$\\
 \hline
\end{tabular}

In the following theorem, we define an auxiliary sequence  $\{A_n\}_{n \in \mathbb{Z}}$ which will be useful in the sequel. Note that this sequence is not a member of the linear space $\mathsf{Fibonacci}^{(k)}$.
The proof of the theorem is
a consequence of 
the Pascal Recursion relation~(\ref{genBinoEqn}).

\begin{theorem}
\label{Bkj:PR}
Let $k \geq 2$
and the sequence $\{A_n\}$
defined as
\begin{equation}
A_n =  \sum_{i \in \mathbb{Z}}{ (-1)^i \left\langle {{n-ik} \choose {i-1}} \right\rangle 2^{n+1-i(k+1)}}
\mbox{ for all } n \in \mathbb{Z}.
\end{equation}
Then $A_0=A_1=A_2=\dots=A_{k-1}=0$,
$A_n = A_{n-1} + A_{n-2} + \dots + A_{n-k}-1$
and
$A_n=2A_{n-1}-A_{n-k-1}$.
\end{theorem}
{\em Proof.}
Note that the above summation in the formula of $A_n$ only has a finite number of non-zero terms.
This is because $\left\langle {{n-ik} \choose {i-1}} \right\rangle = 0$
except for $1 \leq i \leq \frac{n+1}{k+1}$ when $n \geq 0$
and 
$\frac{n+1}{k} \leq i \leq  \frac{n+1}{k+1}$ for $n \leq -1$.
It follows that $A_0=A_1=A_2=\dots=A_{k-1}=0$
and $A_k=-1$.

We have
\begin{align*}
   2 A_{n-1} - A_{n-k-1}
    = & 2  \sum{ (-1)^i \left\langle {{n-1-ik} \choose {i-1}} \right\rangle 2^{n-i(k+1)}}
    \\
    & - \sum{ (-1)^i \left\langle {{n-k-1-ik} \choose {i-1}} \right\rangle 2^{n-k-i(k+1)}}
    \\
     = &  \sum{ (-1)^i \left\langle {{n-1-ik} \choose {i-1}} \right\rangle 2^{n+1-i(k+1)}}
    \\
    &  +  \sum{ (-1)^{i+1} \left\langle {{n-1-(i+1)k} \choose {i-1}} \right\rangle 2^{n+1-(i+1)(k+1)}}
    .
    \end{align*}
  
In the last summation, let $i:=i+1$, we have
\begin{align*}
   2 A_{n-1} - A_{n-k-1}
     = &  \sum{ (-1)^i \left\langle {{n-1-ik} \choose {i-1}} \right\rangle 2^{n+1-i(k+1)}}
    \\
    & +  \sum{ (-1)^i \left\langle {{n-1-ik} \choose {i-2}} \right\rangle 2^{n+1-i(k+1)}}
    \end{align*}
and by the Pascal Recursion~(\ref{genBinoEqn}),
\begin{align*}
   2 A_{n-1} - A_{n-k-1}
    =& \sum{ (-1)^i \left\langle {{n-ik} \choose {i-1}} \right\rangle 2^{n+1-i(k+1)}}
    \\
   =& A_n .
\end{align*}

Therefore, $(\mathtt{R}^{k+1} - 2 \mathtt{R} + \mathtt{I})(A) = 0$.

As $\mathtt{R}^{k+1} - 2 \mathtt{R} + \mathtt{I} = (\mathtt{R} - \mathtt{I})(\mathtt{R}^k + \mathtt{R}^{k-1} + \dots + \mathtt{R} - \mathtt{I})$, it follows that $(\mathtt{R}^k + \mathtt{R}^{k-1} + \dots + \mathtt{R} - \mathtt{I})(A)$ is a constant sequence,
so
$A_{n-1}+A_{n-2}+\dots+A_{n-k}-A_n=
A_{0}+A_{1}+\dots+A_{k-1}-A_k=1$.
\quad
$\blacksquare$

Recall that in Theorem~\ref{BSthm}
we define the sequence 
$S = B^{(0)}+ B^{(1)}+ \dots+ B^{(k-1)}  \in \mathsf{Fibonacci}^{(k)}$.
The following theorem gives an explicit formula for the sequence $S$.

\begin{theorem}
\label{Bkj:theorem1NEW}
Let $k \geq 2$.
The $k$-order Fibonacci sequence $S$ (determined by the first $k$ terms $(1,1,\dots, 1)$) satisfies the following formula
\begin{equation}
\label{formulaSEqn}
S_n = 1 - (k-1) \sum_{i \in \mathbb{Z}}{ (-1)^i \left\langle {{n-ik} \choose {i-1}} \right\rangle 2^{n+1-i(k+1)}}
\mbox{ for all } n \in \mathbb{Z}.
\end{equation}
\end{theorem}
{\em Proof.}
Let $S_n'$ denote the sequence on the RHS of (\ref{formulaSEqn}) then 
$S_n'=1-(k-1)A_n$ where $\{A_n\}$ is the auxiliary sequence defined in
Theorem~\ref{Bkj:PR}.
It follows from Theorem~\ref{Bkj:PR} that
$S_0'=S_1'=\dots=S_{k-1}'=1$, $S_k'=k$ and
$S'_n =2 S'_{n-1} - S'_{n-k-1}$.
By Theorem~\ref{thm_operator}(vi),
the sequence $S$ also satisfies the 
same recursion equation $S_n =2 S_{n-1} - S_{n-k-1}$.
Since
$S_i=S'_i$ for all $0\leq i \leq k$,
it follows that $S_i=S'_i$ for all $i \in \mathbb{Z}$.
\quad
$\blacksquare$

\begin{theorem}
\label{Bkj:theorem1NEW2}
Let $k \geq 2$.
The $k$-order Fibonacci sequence $S$ (determined by the first $k$ terms $(1,1,\dots, 1)$) satisfies the following formula
\begin{equation}
\label{Spositive}
S_n = 1 - (k-1) \sum_{1 \leq i \leq \frac{n+1}{k+1}}{ (-1)^i {{n-ik} \choose {i-1}} 2^{n+1-i(k+1)}}
\mbox{ for all } n \geq 0,
\end{equation}
\begin{equation}
\label{Snegative}
S_n = 1 - (k-1) \sum_{\frac{n+1}{k} \leq i \leq  \frac{n+1}{k+1} }{ (-1)^i \left\langle {{n-ik} \choose {i-1}} \right\rangle 2^{n+1-i(k+1)}}
\mbox{ for all } n \leq -1.
\end{equation}
\end{theorem}
{\em Proof.}
Since $\left\langle {{n-ik} \choose {i-1}} \right\rangle = 0$
except for $1 \leq i \leq \frac{n+1}{k+1}$ when $n \geq 0$
and 
$\frac{n+1}{k} \leq i \leq  \frac{n+1}{k+1}$ for $n \leq -1$,
the theorem follows from Theorem~\ref{Bkj:theorem1NEW}.
\quad
$\blacksquare$

\begin{theorem}
\label{Bkj:theorem2}
Let $k \geq 2$, $0 \leq j \leq k-1$.
The $k$-order Fibonacci sequence $B^{(j)}$ satisfies the following formula
\begin{align*}
B_n^{(j)} = & - \sum_{i \in \mathbb{Z}}{ (-1)^i \left\langle {{n-ik} \choose {i-1}} \right\rangle 2^{n+1-i(k+1)}}
\\
&
+ \sum_{i \in \mathbb{Z}}{ (-1)^i  \left\langle {{n-j-1-ik} \choose {i-1}} \right\rangle 2^{n-j-i(k+1)}}
\mbox{ for all } n \in \mathbb{Z}.
\end{align*}
\end{theorem}
{\em Proof.}
By Theorem~\ref{BSthm}(vii),
$B^{(j)} = \frac{1}{k-1} (\mathtt{I}- \mathtt{R}^{j+1})(S)$,
thus, using the formula~(\ref{formulaSEqn})
for $S_n$ in
Theorem~\ref{Bkj:theorem1NEW},
we obtain 
the desired formula for $B^{(j)}_n$.
    \quad
$\blacksquare$ 

The formula~(\ref{Spositive}) for $S_n$
in Theorem~\ref{Bkj:theorem1NEW2} is equivalent to a formula in Ferguson~\cite{Ferguson1966} (formula (3) for $V_{n,a(n+1)+b}$).
Theorem~\ref{Bkj:theorem2}
for the case $j=k-1$ and positive indices
is proved in Benjamin et al.~\cite{Benjamin2014}.

\section{Explicit formula based on  multinomials}

\label{GF:sec4}

In this section, we will derive 
explicit formula
for 
the two-sided Fibonacci basis sequences
$B^{(0)}, B^{(1)}, \dots, B^{(k-1)}$
expressed in terms of multinomial coefficients.
Since the traditional  multinomial notation is
associated with non-negative integers,
to use these for
our two-sided sequences we need to extend the  multinomial notation
 to include negative integers.
To this end, 
we extend the  multinomial notation ${n \choose {i_1, i_2, \dots, i_t}}$ to negative values of $n$ and $i_1, i_2, \dots, i_t$.

A multinomial is defined as
\begin{align*}
(i_1, i_2, \dots, i_t) & =
    {{i_1 + i_2 + \dots + i_t} \choose {i_1, i_2, \dots, i_t}}
   =
\frac{(i_1 + i_2 + \dots + i_t)!}{i_1! i_2! \dots i_t!}.
\end{align*}

We observe that 
$$
(i_1, i_2, \dots, i_t)=
{{i_1 + \dots + i_t} \choose {i_2 + \dots + i_{t}}}
{{i_2 + \dots + i_{t}} \choose {i_3 + \dots + i_{t}}}
\dots 
{{i_{t-2} + i_{t-1} + i_{t}} \choose {i_{t-1} + i_t}}
{{i_{t-1} + i_{t}} \choose {i_t}}
.
$$

We will use this formula to extend multinomial notation for negative integers.

\begin{definition}
Let $t \geq 2$ be an integer.
For any integers $i_1, i_2, \dots, i_t$,
the generalized  multinomial
$\left\langle (i_1, i_2, \dots, i_t) \right\rangle$ is defined as
\begin{align*}
& \left\langle (i_1, i_2, \dots, i_t) \right\rangle =  \left\langle {{i_1 + i_2 + \dots + i_t} \choose {i_1, i_2, \dots, i_t}} \right\rangle
\\
& =
\left\langle {{i_1 + \dots + i_t} \choose {i_2 + \dots + i_{t}}} \right\rangle
\left\langle {{i_2 + \dots + i_{t}} \choose {i_3 + \dots + i_{t}}} \right\rangle
\dots 
\left\langle {{i_{t-2} + i_{t-1} + i_{t}} \choose {i_{t-1} + i_t}} \right\rangle
\left\langle {{i_{t-1} + i_{t}} \choose {i_t}} \right\rangle
.
\end{align*}
\end{definition}

Using the following formula for 
the generalized binomial coefficient 
\begin{align*}
 \left\langle {{n} \choose {i}} \right\rangle
 & =
  \begin{cases}
  \frac{n^{\underline{n-i}}}{(n-i)!}=
\frac{n(n-1)(n-2)\dots(i+1)}{(n-i)!}, & \text{if } n \geq i\\
            0, & \text{otherwise}
		 \end{cases},
\end{align*}
we obtain the following formula for the 
generalized  multinomial
\begin{align*}
& \left\langle (i_1, i_2, \dots, i_t) \right\rangle =  \left\langle {{i_1 + i_2 + \dots + i_t} \choose {i_1, i_2, \dots, i_t}} \right\rangle
\\
 & =
  \begin{cases}
  \cfrac{(i_1+\dots+i_t)^{\underline{i_1}}
  (i_2+\dots+i_t)^{\underline{i_2}}
  \dots 
  (i_{t-1}+i_t)^{\underline{i_{t-1}}}}{i_1! i_2! \dots i_{t-1}!}, 
  & \text{if } i_1, i_2, \dots, i_{t-1} \geq 0\\
            0, & \text{otherwise}
		 \end{cases}
   .
\end{align*}

When $t=2$,
the Pascal Recursion relation becomes
\begin{align}
\nonumber
& \left\langle (i_1, i_2) \right\rangle =
\left\langle (i_1-1, i_2) \right\rangle
+
\left\langle (i_1, i_2-1) \right\rangle.
\end{align}
For a general $t \geq 2$,
we have the following generalized 
Pascal Recursion relation for multinomials:
\begin{align}
\nonumber
& \left\langle (i_1, i_2, \dots, i_t) \right\rangle 
\\
\label{genBinoEqn2} 
&=
\left\langle (i_1-1, i_2, \dots, i_t) \right\rangle
+
\left\langle (i_1, i_2-1, \dots, i_t) \right\rangle
+ \dots + 
\left\langle (i_1, i_2, \dots, i_t-1) \right\rangle
.
\end{align}

Since $\left\langle {n \choose i} \right\rangle$ is non-zero only for 
$n \geq i \geq 0$ or 
$-1 \geq n \geq i$,
the generalized multinomial
$\left\langle (i_1, i_2, \dots, i_t) \right\rangle$ is non-zero only for
$i_1 + \dots + i_t \geq i_2 + \dots + i_t \geq \dots \geq i_{t-1} + i_t \geq i_t \geq 0$
or $-1 \geq i_1 + \dots + i_t \geq i_2 + \dots + i_t \geq \dots \geq i_{t-1} + i_t \geq i_t$. Using the formula~(\ref{nchoosei2}) for $\left\langle {n \choose i} \right\rangle$,
we can derive the formula for the generalized multinomial in these two separate cases.

Case 1. If $i_1 + \dots + i_t \geq i_2 + \dots + i_t \geq \dots \geq i_{t-1} + i_t \geq i_t \geq 0$, i.e. $i_1, i_2, \dots, i_t \geq 0$, then
\begin{align*}
& \left\langle (i_1, i_2, \dots, i_t) \right\rangle =  \left\langle {{i_1 + i_2 + \dots + i_t} \choose {i_1, i_2, \dots, i_t}} \right\rangle
={{i_1 + i_2 + \dots + i_t} \choose {i_1, i_2, \dots, i_t}}
=(i_1, i_2, \dots, i_t)
.
\end{align*}

Case 2. If $-1 \geq i_1 + \dots + i_t \geq i_2 + \dots + i_t \geq \dots \geq i_{t-1} + i_t \geq i_t$ then
\begin{align*}
 \left\langle (i_1, i_2, \dots, i_t) \right\rangle & =  \left\langle {{i_1 + i_2 + \dots + i_t} \choose {i_1, i_2, \dots, i_t}} \right\rangle
 \\
& =(-1)^{i_1+\dots+i_{t-1}} {{-i_t-1} \choose {i_1, i_2, \dots, i_{t-1}, -i_1-\dots-i_{t}-1}}
\\
&=(-1)^{i_1+\dots+i_{t-1}} (i_1, i_2, \dots, i_{t-1}, -i_1-\dots-i_{t}-1)
.
\end{align*}

Thus, we obtain the following theorem that connects the generalized multinomial
to the classical  multinomial.

\begin{theorem}
For any integer $t \geq 2$ and $i_1,i_2,\dots,i_t \in \mathbb{Z}$, we have
\begin{align*}
& \left\langle (i_1, i_2, \dots, i_t) \right\rangle 
\\
 & =
  \begin{cases}
  (i_1, i_2, \dots, i_t), 
  & \text{if } i_1, i_2, \dots, i_{t} \geq 0\\
  (-1)^{i_1+\dots+i_{t-1}} (i_1, i_2, \dots, i_{t-1}, -i_1-\dots-i_{t}-1)
  & \text{if } i_1, i_2, \dots, i_{t-1} \geq 0 \mbox{ and } i_1 + \dots + i_t \leq -1
  \\
            0, & \text{otherwise.}
		 \end{cases}
\end{align*}
\end{theorem}

In the following theorem, we define an auxiliary sequence $\{X_n\}_{n \in \mathbb{Z}}$. Note that $X$ is a member of the linear space $ \mathsf{Fibonacci}^{(k)}$.

\begin{theorem}
\label{MilesExtendExtra}
Let $k \geq 2$, $c \in \mathbb{Z}$ any constant,
and
\begin{align*}
X_n 
&=
\sum_{a_1 + 2 a_2 + \dots + k a_k = n + c}{\left\langle (a_1, a_2, \dots, a_k) \right\rangle}
\\
&=
\sum_{s_1 + s_2 + \dots + s_k = n + c}{
\left\langle {s_1 \choose s_2} \right\rangle
\left\langle {s_2 \choose s_3} \right\rangle
\dots
\left\langle {s_{k-1} \choose s_k} \right\rangle
}
.
\end{align*}
Then $\{X_n\}_{n \in \mathbb{Z}}$ is a Fibonacci sequence of order $k$.
\end{theorem}
{\em Proof.}
The two formulas on the RHS are equivalent by using the variables $s_1=a_1+\dots+a_k$,
$s_2=a_2+\dots+a_k$, \dots, $s_{k-1}=a_{k-1}+a_k$ and $s_k=a_k$.

Note that the summation only has a finite number of non-zero terms.
This is because $\left\langle (a_1, a_2, \dots, a_k) \right\rangle$ is non-zero
only if $s_1 \geq s_2 \geq \dots \geq s_k \geq 0$ or $-1 \geq s_1 \geq s_2 \geq \dots \geq s_k$, and there are only a finite number 
of choices for $s_1, s_2, \dots, s_k$
that have the same sign
whose sum $s_1 + s_2 + \dots + s_k=n+c$ is fixed.

By Pascal Recursion relation~(\ref{genBinoEqn2}),
\begin{align*}
    X_n =
& 
\sum_{a_1 + 2 a_2 + \dots + k a_k = n + c}{\left\langle (a_1-1, a_2, \dots, a_k) \right\rangle}
\\
& + 
\sum_{a_1 + 2 a_2 + \dots + k a_k = n + c}{\left\langle (a_1, a_2-1, \dots, a_k) \right\rangle}
\\
& + \dots + 
\sum_{a_1 + 2 a_2 + \dots + k a_k = n + c}{\left\langle (a_1, a_2, \dots, a_k-1) \right\rangle}
.
\end{align*}

Let $a_1'=a_1-1$, $a_2'=a_2-1$, \dots, $a_k'=a_k-1$. We have
\begin{align*}
X_n= &
\sum_{a_1' + 2 a_2 + \dots + k a_k = n + c - 1}{\left\langle (a_1', a_2, \dots, a_k) \right\rangle}
\\
& + 
\sum_{a_1 + 2 a_2' + \dots + k a_k = n + c -2}{\left\langle (a_1, a_2', \dots, a_k) \right\rangle}
\\
& + \dots + 
\sum_{a_1 + 2 a_2 + \dots + k a_k' = n + c - k}{\left\langle (a_1, a_2, \dots, a_k') \right\rangle}
\\
& = X_{n-1} + X_{n-2} + \dots + X_{n-k},
\end{align*}
therefore,
$\{X_n\}$ is a Fibonacci sequence of order $k$.
    \quad
$\blacksquare$ 

\begin{theorem}
\label{MilesExtend3}
Let $k \geq 2$. Then
\begin{equation}
\label{eqn_MilesExtend3}
B_{n}^{(0)}
=
\sum_{a_1 + 2 a_2 + \dots +  ka_k = n-k}{\left\langle (a_1, a_2, \dots, a_k) \right\rangle}, \mbox{ for all } n \in \mathbb{Z}.
\end{equation}
\end{theorem}
{\em Proof.}
Let $B'$ denote the RHS, then 
by Theorem~\ref{MilesExtendExtra},
$B'$ is a Fibonacci sequence.
We only need to show its initial values
match  with those of $B^{(0)}$.

Again, as in the proof of Theorem~\ref{MilesExtendExtra}, we use the variables $s_1=a_1+\dots+a_k$,
$s_2=a_2+\dots+a_k$, \dots, $s_{k-1}=a_{k-1}+a_k$ and $s_k=a_k$,
then $s_1+s_2+\dots+s_k=n-k$.
When $n=0$, $s_1 + s_2 + \dots + s_k=-k <0$, so $\left\langle (a_1, a_2, \dots, a_k) \right\rangle$ is non-zero
only if $-1 \geq s_1 \geq s_2 \geq \dots \geq s_{k}$. The only possibility is $s_1=s_2=\dots=s_k=-1$
and this gives $a_1=a_2=\dots=a_{k-1}=0$, $a_k=-1$
and  $B'_0 = \left\langle (0, \dots, 0, -1) \right\rangle = 1$.

When $1 \leq n \leq k-1$,
$-(k-1) \leq s_1 + s_2 + \dots + s_k=n-k <0$.
There are no such $-1 \geq s_1 \geq s_2 \geq \dots \geq s_{k}$ that satisfy this condition, so the summation is empty
and $B'_n = 0$ for $1 \leq n \leq k-1$. 
 \quad
$\blacksquare$

\begin{theorem}
\label{MilesExtend2}
Let $k \geq 2$. Then
$$
B_{n}^{(k-1)}
=
\sum_{a_1 + 2 a_2 + \dots + k a_k = n-k+1}{\left\langle (a_1, a_2, \dots, a_k) \right\rangle}, \mbox{ for all } n \in \mathbb{Z}
.
$$
\end{theorem}
{\em Proof.}
By Theorem~\ref{BSthm}(iii),
$B^{(k-1)}= \mathtt{L}(B^{(0)})$,
so using the formula for $B^{(0)}_n$ in
Theorem~\ref{MilesExtend3}
we obtain the desired formula for $B^{(k-1)}_n$.
\quad
$\blacksquare$ 

The formula in
Theorem~\ref{MilesExtend2}
is proved in Miles~\cite{Miles1960}
for natural number $n \geq k-1$.
Our Theorem~\ref{MilesExtend2} extends it to $n<k-1$ and negative integer $n$.

The Tribonacci sequence $\{T_n \}_{n \geq 0}$
studied in Rabinowitz~\cite{Rabinowitz1996}
is a Fibonacci sequence of order $k=3$
with initial values 
$T_0=0$, $T_1=1$, $T_2=1$.
Solving for $T_{-1}$, we have
$T_{-1} = 0$, so $T = \mathtt{L}(B^{(2)})$.
The formula in
Theorem~\ref{MilesExtend2}
is proved in Rabinowitz~\cite{Rabinowitz1996}
for $k=3$ and $n \geq 2$.
Our Theorem~\ref{MilesExtend2} extends it to all order $k \geq 2$ and all index $n \in \mathbb{Z}$.

The next theorem give 
an explicit formula for all basis Fibonacci sequences of order $k$.

\begin{theorem}
\label{MilesExtend}
Let $k \geq 2$.
For any $0 \leq j \leq k-1$,
$$
B_{n}^{(j)}
=
\sum_{n-k-j \leq  a_1 + 2 a_2 + \dots + k a_k \leq n-k}{\left\langle (a_1, a_2, \dots, a_k) \right\rangle}, \mbox{ for all } n \in \mathbb{Z}
.
$$
\end{theorem}
{\em Proof.}
By Theorem~\ref{BSthm}(ii),
 $B^{(j)} = 
\sum_{i=0}^{j} \mathtt{R}^i(B^{(0)})
$,
so using the formula for $B^{(0)}_n$ in
Theorem~\ref{MilesExtend3}
we obtain the desired formula for $B^{(j)}_n$.
\quad
$\blacksquare$ 

Theorem~\ref{MilesExtend2}
and Theorem~\ref{MilesExtend}
give rise to two different formulas 
for the sequence $B^{(k-1)}$.
It would be interesting to see a combinatorial
proof of the equality of these two formulas.

\section{A remark on a tiling problem}

\label{GF:sec5}

It is well known that the classical Fibonacci sequence, $F_0=0$, $F_1=1$, $F_n=F_{n-1}+F_{n-2}$, has a close relation with the tiling problem.
The value $F_n$ counts the number of
tilings of an $1 \times n$-board with square-tiles $1 \times 1$ and domino-tiles $1 \times 2$.
This is because for $n \geq 2$,
by considering the first tile,
if the first tile is a square then
there are $F_{n-1}$ ways to cover 
the remaining strip of length $n-1$,
and if the first tile is a domino then there are $F_{n-2}$ ways to cover 
the remaining strip of length $n-2$.
That is how the recursion equation
$F_n=F_{n-1}+F_{n-2}$ arises.

If we allow tiles of length up to $k$, then the result is a sequence $\{C_n\}_{n \geq 0}$. We have
$C_0=0$, $C_1=1$,
$C_2=C_0+C_1$,
$C_3=C_0+C_1+C_2$,\dots,
$C_{k-1}=C_0+C_1+\dots +C_{k-2}$,
and for $n \geq k$,
$C_{n}=C_{n-1}+C_{n-2}+\dots +C_{n-k}$.
Of course, if we extend the index to negative integers and set 
$C_{-1}=C_{-2}=\dots=C_{-(k-2)}=0$
then we have the Fibonacci recursion equation $C_{n}=C_{n-1}+C_{n-2}+\dots +C_{n-k}$ holds for all $n \geq 2$.
This sequence $C$ is just a left shift
of the basis sequence $B^{(k-1)}$.
Indeed, $C = \mathtt{R}^{k-2}(B^{(k-1)})$.
Many authors such as Gabai, Philippou, Muwafi,
Benjamin, Heberle, Quinn and Su
~\cite{Philippou1982,Gabai1970,Benjamin2014,Benjamin2000} have studied this
tiling problem and here we decide to use the letter $C$ to denote this sequence
since it is related to a {\em combinatorial} problem.

\end{document}